\documentclass[final,1p,times,11.5pt]{article}

\usepackage{enumerate}
\usepackage{epsfig,alltt}

\usepackage{amssymb}
\usepackage{amsfonts}
\usepackage{graphicx}
\usepackage{amsmath}
\usepackage{amsthm}

\usepackage[utf8]{inputenc}
\usepackage[T1]{fontenc}

\usepackage{floatflt}
\usepackage{makeidx}
\usepackage{algorithm}

\usepackage{algorithmic}

\usepackage{url}
\usepackage{float} 

\theoremstyle{remark}

\usepackage[english]{babel}

\usepackage{amsfonts}
\usepackage{amssymb}
\usepackage{geometry}
\usepackage{color}

\newcommand{\begeq}[1]{\begin{equation} \label{#1}}
\newcommand{\fineq}{\end{equation}}

\newcommand{\g}{\gamma}

\def\bR{{\Bbb R}}                  

\title{Large and moderate deviation principles for recursive kernel density estimators defined by stochastic approximation method}
\author{Yousri Slaoui}

\begin{document}

\newtheorem{theor}{Theorem}
\newtheorem{prop}{Proposition}
\newtheorem{lemma}{Lemma}
\newtheorem{lem}{Lemma}
\newtheorem{coro}{Corollary}

\newtheorem{prof}{Proof}
\newtheorem{defi}{Definition}
\newtheorem{rem}{Remark}

\date{ }
\maketitle

\begin{abstract}
In this paper we prove large and moderate deviations principles
for the recursive kernel estimators of a probability density function defined by the stochastic approximation algorithm introduced by Mokkadem et al. [2009. The stochastic approximation method for the estimation of a probability density. J. Statist. Plann. Inference 139, 2459-2478]. We show that the estimator constructed using the stepsize which minimize the variance of the class of the recursive estimators defined in Mokkadem et al. (2009) gives the same pointwise LDP and MDP as the Rosenblatt kernel estimator. We provide results both for the pointwise and the uniform deviations.
\end{abstract}

2000 \textit{Mathematics Subject Classifiation}: $62$G$07$, $62$L$20$  , $60$F$10$.\\
\textbf{Key Words:} Density estimation; Stochastic approximation algorithm; Large and Moderate deviations principles.

\section{Introduction}

Let $X_1,\ldots ,X_n$ be independent, identically distributed 
$\bR^d$-valued random vectors, 
and let $f$ denote the probability density of $X_1$. To construct a stochastic 
algorithm, which approximates the function $f$ at a given point $x$, Mokkadem et al. (2009) defined an algorithm of search of the zero of the function $h~: y\mapsto f(x)-y$. They 
proceed as follows: (i) they set $f_0(x)\in\bR$; 
(ii) for all $n\geq 1$, they set  
\begin{eqnarray*}
f_n(x)=f_{n-1}(x)+\g_nW_n(x)
\end{eqnarray*}
where $W_n(x)$ is an ``observation'' of the function $h$ at the point 
$f_{n-1}(x)$ and $(\g_n)$ is a sequence of positive real numbers that goes to zero. To define $W_n(x)$, they follow the approach of R\'ev\'esz (1973, 1977) and of Tsybakov (1990), and introduced a kernel $K$ (which is a function satisfying  
$\int_{\bR^d} K(x)dx=1$) and a bandwidth $(h_n)$ (which is a sequence of positive real 
numbers that goes to zero), and they set
$W_n(x)=h_n^{-d}K(h_n^{-1}[x-X_n])-f_{n-1}(x)$. 
The stochastic approximation algorithm introduced in Mokkadem et al. (2009) which estimate recursively the density $f$ at the point $x$ is
\begin{eqnarray}\label{algo:density}
f_n(x)=(1-\g_n)f_{n-1}(x)+\g_nh_n^{-d}K\left(\frac{x-X_n}{h_n}\right).
\end{eqnarray}
Recently, large and moderate deviations results have been proved for the
well-known nonrecursive kernel density estimator introduced by Rosenblatt (1956) (see also Parzen, 1962). The large deviations principle
has been studied by Louani (1998) and Worms (2001). Gao (2003) and Mokkadem et al. (2005) extend these results and provide moderate deviations principles. The purpose of this paper is to establish large and moderate deviations principles for the recursive density estimator defined by the stochastic approximation algorithm~(\ref{algo:density}).

Let us first recall that a $\mathbb{R}^m$-valued sequence $\left(Z_n\right)_{n\geq 1}$ satisfies a large deviations principle (LDP) with speed $\left(\nu_n\right)$ and good rate function $I$ if :

\begin{enumerate}
\item $\left(\nu_n\right)$ is a positive sequence such that $\lim_{n\to \infty}\nu_n=\infty$;
\item $I:\mathbb{R}^m\to \left[0,\infty\right]$ has compact level sets;
\item for every borel set $B\subset \mathbb{R}^m$,
\begin{eqnarray*}
-\inf_{x\in \overset{\circ}{B}}I\left(x\right)&\leq & \liminf_{n\to \infty} \nu_n^{-1}\log \mathbb{P}\left[Z_n\in B\right]\\
&\leq & \limsup_{n\to \infty} \nu_n^{-1}\log \mathbb{P}\left[Z_n\in B\right]\leq -\inf_{x\in \overline{B}}I\left(x\right),
\end{eqnarray*}
 where $\overset{\circ}{B}$ and $\overline{B}$ denote the interior and the closure of $B$ respectively. Moreover, let $\left(v_n\right)$ be a nonrandom sequence that goes to infinity; if $\left(v_nZ_n\right)$ satisfies a LDP,
then $\left(Z_n\right)$ is said to satisfy a moderate deviations principle (MDP).
\end{enumerate}

The first aim of this paper is to establish pointwise LDP for the recursive kernel density estimators defined by the stochastic approximation algorithm~(\ref{algo:density}). It turns out that the rate function depend on the choice of the stepsize $\left(\gamma_n\right)$; In the first part of this paper we focus on the following two special cases : (1) $\left(\gamma_n\right)=\left(n^{-1}\right)$ and (2) $\left(\gamma_n\right)=\left(h_n^d\left(\sum_{k=1}^nh_k^d\right)^{-1}\right)$, the first one belongs to the subclass of recursive kernel estimators which have a minimum MSE or MISE and the second choice belongs to the subclass of recursive kernel estimators which have a minimum variance (see Mokkadem et al., 2009).\\

We show that using the stepsize $\left(\gamma_n\right)=\left(n^{-1}\right)$ and $\left(h_n\right)\equiv \left(cn^{-a}\right)$ with $c>0$ and $a\in \left]0,1/d\right[$, the sequence $\left(f_n\left(x\right)-f\left(x\right)\right)$ satisfies a LDP with speed $\left(nh_n^d\right)$ and the rate function defined as follows:
\begin{eqnarray}\label{sys:LDPWW}
\left\{\begin{array}{llll}
if & f\left(x\right)\not=0, &I_{a,x}:t\to f\left(x\right)I_a\left(1+\frac{t}{f\left(x\right)}\right)\\
if & f\left(x\right)=0, &I_{a,x}\left(0\right)=0\quad and \quad I_{a,x}\left(t\right)=+\infty \quad for\quad t\not=0.
\end{array}
\right.
\end{eqnarray}
where 
\begin{eqnarray*}
&&I_a\left(t\right)=\sup_{u\in \mathbb{R}}\left\{ut-\psi_a\left(u\right)\right\}\\
&&\psi_a\left(u\right)=\int_{\left[0,1\right]\times\mathbb{R}^d}s^{-ad}\left(e^{us^{ad}K\left(z\right)}-1\right)dsdz,
\end{eqnarray*}
which is the same rate function for the LDP of the Wolverton and Wagner (1969) kernel estimator (see Mokkadem et al., 2006).\\

Moreover, we show that using the stepsize $\left(\gamma_n\right)=\left(h_n^d\left(\sum_{k=1}^nh_k^d\right)^{-1}\right)$ and more general bandwiths defined as $h_n=h\left(n\right)$ for all $n$, where $h$ is a regularly varing function with exponent $\left(-a\right)$, $a\in \left]0,1/d\right[$. We prove that the sequence $\left(f_n\left(x\right)-f\left(x\right)\right)$ satisfies a LDP with speed $\left(nh_n^d\right)$ and the rate function defined as follows:
\begin{eqnarray}\label{sys:LDPDe}
\left\{\begin{array}{llll}
if & f\left(x\right)\not=0, &I_{x}:t\to f\left(x\right)I\left(1+\frac{t}{f\left(x\right)}\right)\\
if & f\left(x\right)=0, &I_{x}\left(0\right)=0\quad and \quad I_{x}\left(t\right)=+\infty \quad for\quad t\not=0.
\end{array}
\right.
\end{eqnarray}
where
\begin{eqnarray*}
&&I\left(t\right)=\sup_{u\in \mathbb{R}}\left\{ut-\psi\left(u\right)\right\}\\
&&\psi\left(u\right)=\int_{\mathbb{R}^d}\left(e^{uK\left(z\right)}-1\right)dz,
\end{eqnarray*}
which is the same rate function for the LDP of the Rosenblatt kernel estimator (see Mokkadem et al., 2005).\\

Our second aim is to provide pointwise MDP for the density estimator defined by the stochastic approximation algorithm~(\ref{algo:density}). In this case, we consider more general stepsizes defined as $\gamma_n=\gamma\left(n\right)$ for all $n$, where $\gamma$ is a regularly function with exponent $\left(-\alpha\right)$, $\alpha\in \left]1/2,1\right]$.
Throughout this paper we will use the following notation:
\begin{eqnarray}
&&\xi= \lim_{n\to +\infty}\left(n\gamma_n\right)^{-1}.\label{not:xi}
\end{eqnarray}
For any positive sequence $\left(v_n\right)$ satisfying 
\begin{eqnarray*}
\lim_{n\to \infty}v_n=\infty \quad \mbox{and} \quad \lim_{n\to \infty}\frac{\gamma_nv_n^2}{h_n^{d}}=0
\end{eqnarray*}
and general bandwidths $\left(h_n\right)$, we prove that the sequence
\begin{eqnarray*}
v_n\left(f_n\left(x\right)-f\left(x\right)\right)
\end{eqnarray*} 
satisfies a LDP of speed $\left(h_n^d/\left(\gamma_nv_n^2\right)\right)$ and rate function $J_{a,\alpha,x}\left(.\right)$ defined by 
\begin{eqnarray}\label{sys:MDP}
\left\{\begin{array}{llll}
if & f\left(x\right)\not=0, &J_{a,\alpha,x}:t\to \frac{t^2\left(2-\left(\alpha-ad\right)\xi\right)}{2f\left(x\right)\int_{\mathbb{R}^d}K^2\left(z\right)dz}\\
if & f\left(x\right)=0, &J_{a,\alpha,x}\left(0\right)=0\quad and \quad J_{a,\alpha,x}\left(t\right)=+\infty \quad for\quad t\not=0.
\end{array}
\right.
\end{eqnarray}
Let us point out that using the stepsize $\left(\gamma_n\right)=\left(h_n^d\left(\sum_{k=1}^dh_k^d\right)^{-1}\right)$ which minimize the variance of $f_n$, we obtain the same rate function for the pointwise LDP and MDP as the one obtained for the Rosenblatt kernel estimator.    \\
Finally, we give a uniform version of the previous results. More precisely,
let $U$ be a subset of $\mathbb{R}^d$; we establish large and moderate deviations principles for the sequence $\left(\sup_{x\in U}\left|f_n\left(x\right)-f\left(x\right)\right|\right)$. 

\section{Assumptions and main results} \label{section 2}
We define the following class of regularly varying sequences.
\begin{defi}
Let $\gamma \in \mathbb{R} $ and $\left(v_n\right)_{n\geq 1}$ be a nonrandom positive sequence. We say that $\left(v_n\right) \in \mathcal{GS}\left(\gamma \right)$ if
\begin{eqnarray}\label{eq:5}
\lim_{n \to +\infty} n\left[1-\frac{v_{n-1}}{v_{n}}\right]=\gamma .
\end{eqnarray}
\end{defi}
Condition~(\ref{eq:5}) was introduced by Galambos and Seneta (1973) to define regularly varying sequences (see also Bojanic and Seneta, 1973), and by Mokkadem and Pelletier (2007) in the context of stochastic approximation algorithms. Typical sequences in $\mathcal{GS}\left(\gamma \right)$ are, for $b\in \mathbb{R}$, $n^{\gamma}\left(\log n\right)^{b}$, $n^{\gamma}\left(\log \log n\right)^{b}$, and so on. \\

\subsection{Pointwise LDP for the density estimator defined by the stochastic approximation algorithm~(\ref{algo:density})}
\subsubsection{Choices of $\left(\gamma_n\right)$ minimizing the MISE of $f_n$} 
It was shown in Mokkadem et al. (2009) that to minimize the MISE of $f_n$, the stepsize $(\gamma_n)$ must be chosen in $\mathcal{GS}\left(-1\right)$ and must satisfy $\lim_{n\to \infty} n\gamma_n= 1$. The most simple example of stepsize belonging to $\mathcal{GS}\left(-1\right)$ and such that $\lim_{n\to \infty} n\gamma_n= 1$ is $\left(\gamma_n\right)=\left(n^{-1}\right)$. For this choice of
stepsize, the estimator $f_n$ defined by~(\ref{algo:density}) equals the recursive kernel estimator introduced by Wolverton and Wagner (1969).\\

To establish pointwise LDP for $f_n$ in this case, we need the following assumptions.
\begin{description}
\item(L1) $K:\mathbb{R}^d\rightarrow \mathbb{R}$ is a bounded and integrable function satisfying $\int_{\mathbb{R}^d}K\left( z\right) dz=1$, and $\lim_{\|z\|\to \infty}K\left(z\right)=0$. 
\item(L2) i) $\left(h_n\right)=\left(cn^{-a}\right)$ with $a\in \left]0,1/d\right[$ and $c>0$. \\
  $ii)$ $\left(\gamma_{n}\right)=\left(n^{-1}\right)$.
\end{description}
The following Theorem gives the pointwise LDP for $f_n$ in this case.
\begin{theor}[Pointwise LDP for Wolverton and Wagner estimator]\label{the:LDPWW} $ $\\
Let Assumptions $\left(L1\right)$ and $\left(L2\right)$ hold and assume that $f$ is continuous at $x$. Then, the sequence $\left(f_n\left(x\right)-f\left(x\right)\right)$ satisfies a LDP with speed $\left(nh_n^d\right)$ and rate function defined by~(\ref{sys:LDPWW}).
\end{theor}

\subsubsection{Choices of $\left(\gamma_n\right)$ minimizing the variance of $f_n$}
It was shown in Mokkadem et al. (2009) that to minimize the asymptotic variance of $f_n$, the stepsize $(\gamma_n)$ must be chosen in $\mathcal{GS}\left(-1\right)$ and must satisfy $\lim_{n\to \infty} n\gamma_n= 1-ad$. The most simple example of stepsize belonging to $\mathcal{GS}\left(-1\right)$ and such that $\lim_{n\to \infty} n\gamma_n= 1-ad$ is $\left(\gamma_n\right)=\left(\left(1-ad\right)n^{-1}\right)$, an other stepsize satisfying this conditions is $\left(\gamma_n\right)=\left(h_n^d\left(\sum_{k=1}^nh_k^d\right)^{-1}\right)$. For this last choice of stepsize, the estimator $f_n$ defined by~(\ref{algo:density}) produces the estimator considered by Deheuvels (1973) and Duflo (1997).\\
To establish pointwise LDP for $f_n$ in this case, we assume that.
\begin{description}
\item(L3) i) $\left(h_{n}\right)\in \mathcal{GS} \left(-a\right)$ with 
$a \in \left]0,1/{d}\right[$.\\
  $ii)$ $\left(\gamma_{n}\right)=\left(h_n^d\left(\sum_{k=1}^nh_k^d\right)^{-1}\right)$.
\end{description}
The following Theorem gives the pointwise LDP for $f_n$ in this case.
\begin{theor}[Pointwise LDP for Deheuvels estimator]\label{the:LDPDeu} $ $\\
Let Assumptions $\left(L1\right)$ and $\left(L3\right) $ hold and assume that $f$ is continuous at $x$. Then, the sequence $\left(f_n\left(x\right)-f\left(x\right)\right)$ satisfies a LDP with speed $\left(nh_n^d\right)$ and rate function defined by~(\ref{sys:LDPDe}).

\end{theor}

\subsection{Pointwise MDP for the density estimator defined by the stochastic approximation algorithm~(\ref{algo:density})}
Let $\left(v_n\right)$ be a positive sequence; we assume that
\begin{description} 
\item(M1) $K:{\mathbb R^d}\rightarrow {\mathbb R}$ is a continuous, bounded 
function satisfying $\int_{\mathbb{R}^d}K\left( z\right) dz=1$, and, for 
all $j\in\{1,\ldots d\}$, $\int_{\mathbb{R}}z_jK\left( z\right) dz_j=0$ 
and $\int_{\mathbb{R}^d}z_j^2|K\left( z\right)| dz<\infty$. 
\item(M2) $i)$  $\left(\gamma_n\right) \in \mathcal{GS}\left(-\alpha \right)$ with 
$\alpha\in \left]{1}/{2},1\right]$.  \\
  $ii)$ $\left(h_{n}\right)\in \mathcal{GS} \left(-a\right)$ with 
$a \in \left]0,{\alpha}/{d}\right[$.\\ 
 $iii)$  $\lim_{n\to\infty}\left(n\gamma_n\right)\in]\min\{2a,(\alpha-ad)/2\},\infty]$.
\item(M3) $f$ is bounded, twice differentiable, and, for all $i,j\in\{1,\ldots d\}$, 
$\partial^2f/\partial x_i\partial x_j$ is bounded.
\item(M4) $\lim_{n\to \infty}v_n=\infty$ and $\lim_{n\to \infty}\gamma_nv_n^2/h_n^d=0$.
\end{description}
The following Theorem gives the pointwise MDP for $f_n$.
\begin{theor}[Pointwise MDP for the recursive estimators defined by~(\ref{algo:density})]\label{the:MDP} $ $\\
Let Assumptions $\left(M1\right)-\left(M4\right) $ hold and assume that $f$ is continuous at $x$. Then, the sequence $\left(f_n\left(x\right)-f\left(x\right)\right)$ satisfies a MDP with speed $\left(h_n^d/\left(\gamma_nv_n^2\right)\right)$ and rate function $J_{a,\alpha,x}$ defined in (\ref{sys:MDP}).
\end{theor}
\subsection{Uniform LDP and MDP for the density estimator defined by the stochastic approximation algorithm~(\ref{algo:density})}
To establish uniform large deviations principles for the density
estimator defined by the stochastic approximation algorithm~(\ref{algo:density}) on a bounded set, we need the following assumptions:
\begin{description} 
\item(U1) $i)$ For all $j\in\{1,\ldots d\}$, $\int_{\mathbb{R}}z_jK\left( z\right) dz_j=0$ 
and $\int_{\mathbb{R}^d}z_j^2|K\left( z\right)| dz<\infty$.\\
$ii)$ $K$ is H\"older continuous.
\item(U2) $f$ is bounded, twice differentiable, and, $\sup_{x\in\mathbb{R}^d}\|D^2f\left(x\right)\|<\infty$.\\
\item(U3) $\lim_{n\to\infty}\frac{\gamma_nv_n^2\log\left(1/h_n\right)}{h_n^d}=0$ and $\lim_{n\to\infty}\frac{\gamma_nv_n^2\log v_n}{h_n^d}=0$.
\end{description}

Set $U\subseteq \mathbb{R}^{d}$; in order to state in a compact form the uniform large and moderate deviations principles for the density estimator defined by the stochastic approximation algorithm~(\ref{algo:density}) on $U$,
we set:
\begin{eqnarray*}
g_U\left(\delta\right)&=&\left\{\begin{array}{lllll}
\|f\|_{U,\infty}I_a\left(1+\frac{\delta}{\|f\|_{U,\infty}}\right) & \mbox{when} & v_n\equiv 1& \mbox{,}&\left(L1\right)\,\,\mbox{and}\,\,\left(L2\right)\,\,\mbox{hold}\\
\|f\|_{U,\infty}I\left(1+\frac{\delta}{\|f\|_{U,\infty}}\right) & \mbox{when} & v_n\equiv 1&\mbox{,}&\left(L1\right)\,\,\mbox{and}\,\,\left(L3\right)\,\,\mbox{hold}\\
\frac{\delta^2\left(2-\left(\alpha-ad\right)\xi\right)}{2\|f\|_{U,\infty}\int_{\mathbb{R}^d}K^2\left(z\right)dz}
 & \mbox{when} & v_n\to \infty&\mbox{,}& \left(M1\right)-\left(M4\right)\,\,\mbox{hold} \\
\end{array}
\right.\\
\tilde{g}_U\left(\delta\right)&=&\min\left\{g_U\left(\delta\right),g_U\left(-\delta\right)\right\}
\end{eqnarray*}
where $\|f\|_{U,\infty}=\sup_{x\in U}\left|f\left(x\right)\right|$.
\begin{rem}
The functions $g_{U}\left(.\right)$ and $\tilde{g}_{U}\left(.\right)$ are non-negative, continuous, increasing on $\left]0,+\infty\right[$ and decreasing on $\left]-\infty,0\right[$, with a unique global minimum in $0$ ($\tilde{g}_U\left(0\right)=g_U\left(0\right)=0$). They are thus good rate functions (and $g_U(.)$ is strictly convex).
\end{rem}
Theorem~\ref{theo:unifbounded} below states uniform LDP on $U$ in the case $U$ is bounded, and Theorem~\ref{theo:unifunbounded} in the case $U$ is unbounded.
\begin{theor}[Uniform deviations on a bounded set for the recursive estimator defined by~(\ref{algo:density})]\label{theo:unifbounded}
Let $\left(U1\right)-\left(U3\right)$ hold. Then for any bounded subset $U$ of $\mathbb{R}^d$ and for all $\delta>0$, 
\begin{eqnarray}\label{eq:unifboun}
\lim_{n\to \infty} \frac{\gamma_nv_n^2}{h_n^d}\log \mathbb{P}\left[\sup_{x\in U}v_n\left|f_n\left(x\right)-f\left(x\right)\right|\geq \delta\right]=-\tilde{g}_U\left(\delta\right)
\end{eqnarray} 
\end{theor}
To establish uniform large deviations principles for the density
estimator defined by the stochastic approximation algorithm~(\ref{algo:density}) on an unbounded set, we need the following additionnal assumptions:
\begin{description} 
\item(U4) $i)$ There exists $\beta >0$ such that $\int_{\mathbb{R}^{d}}\|x\|^{\beta}f\left(x\right)dx<\infty$.\\
$ii)$ $f$ is uniformly continuous.
\item(U5) There exists $\tau >0$ such that $z\mapsto \|z\|^{\tau}K\left(z\right)$ is a bounded function.
\item(U6) $i)$ There exists $\zeta>0$ such that $\int_{\mathbb{R}^d}\|z\|^{\zeta}\left|K\left(z\right)\right|dz<\infty$\\
$ii)$ There exists $\eta>0$ such that $z\mapsto \|z\|^{\eta}f\left(z\right)$ is a bounded function.
\end{description}
\begin{theor}[Uniform deviations on an unbounded set for the recursive estimator defined by~(\ref{algo:density})]\label{theo:unifunbounded}
Let $\left(U1\right)-\left(U6\right)$ hold. Then for any subset $U$ of $\mathbb{R}^d$ and for all $\delta>0$, 
\begin{eqnarray*}
-\tilde{g}_U\left(\delta\right)&\leq &\liminf_{n\to \infty} \frac{\gamma_nv_n^2}{h_n^d}\log \mathbb{P}\left[\sup_{x\in U}v_n\left|f_n\left(x\right)-f\left(x\right)\right|\geq \delta\right]\\
&\leq&\limsup_{n\to \infty} \frac{\gamma_nv_n^2}{h_n^d}\log \mathbb{P}\left[\sup_{x\in U}v_n\left|f_n\left(x\right)-f\left(x\right)\right|\geq \delta\right]\leq -\frac{\beta}{\beta+d}\tilde{g}_U\left(\delta\right)
\end{eqnarray*} 
\end{theor}
The following corollary is a straightforward consequence of Theorem~\ref{theo:unifunbounded}.
\begin{coro}\label{coro:unifunbounded}
Under the assumptions of Theorem~\ref{theo:unifunbounded}, if $\int_{\mathbb{R}^d}\|x\|^{\xi}f\left(x\right)dx<\infty$ for all $\xi$ in $\mathbb{R}$, then for any subset $U$ of $\mathbb{R}^d$,
\begin{eqnarray}\label{eq:unifunboun}
\lim_{n\to \infty} \frac{\gamma_nv_n^2}{h_n^d}\log \mathbb{P}\left[\sup_{x\in U}v_n\left|f_n\left(x\right)-f\left(x\right)\right|\geq \delta\right]=-\tilde{g}_U\left(\delta\right)
\end{eqnarray}
\end{coro}
\paragraph{Comment.} Since the sequence $\left(\sup_{x\in U}\left|f_n\left(x\right)-f\left(x\right)\right|\right)$ is positive and since $\tilde{g}_U$ is continuous on $\left[0,+\infty\right[$, increasing and goes to infinity as $\delta\to \infty$, the application of Lemma 5 in Worms (2001) allows to deduce from~(\ref{eq:unifboun}) or~(\ref{eq:unifunboun}) that $\sup_{x\in U}\left|f_n\left(x\right)-f\left(x\right)\right|$ satisfies a LDP with speed $\left(\gamma_n^{-1}h_n^d\right)$ and good rate function $\tilde{g}_U$ on $\mathbb{R}^+$.
\section{Proofs}
Throught this section we use the following notation:
\begin{eqnarray}
&&\Pi_n=\prod_{j=1}^n\left(1-\gamma_j\right),\nonumber\\
&&Z_n\left(x\right)=h_n^{-d}Y_n,\nonumber\\
&&Y_n=K\left(\frac{x-X_n}{h_n}\right)\label{eq:Y}
\end{eqnarray}
Throughout the proofs, we repeatedly apply Lemma 2 in Mokkadem et al. (2009). For the convenience of the reader, we state it now.
\begin{lemma}\label{lemma:Tech} Let $\left(v_n\right)\in \mathcal{GS}\left(v^{*}\right)$, $\left(\gamma_n\right)\in \mathcal{GS}\left(-\alpha\right)$, and $m>0$ such that $m-v^*\xi>0$ where $\xi$ is defined in~(\ref{not:xi}). We have
\begin{eqnarray*}
\lim_{n\to +\infty}v_n\Pi_n^m\sum_{k=1}^n\Pi_k^{-m}\frac{\gamma_k}{v_k}=\frac{1}{m-v^*\xi}.
\end{eqnarray*} 
Moreover, for all positive sequence $\left(\alpha_n\right)$ such that $\lim_{n\to +\infty}\alpha_n=0$, and for all $\delta\in \mathbb{R}$,
\begin{eqnarray*}
\lim_{n\to +\infty}v_n\Pi_n^m\left[\sum_{k=1}^n\Pi_k^{-m}\frac{\gamma_k}{v_k}\alpha_k+\delta\right]=0.
\end{eqnarray*}  
\end{lemma}
Noting that, in view of~(\ref{algo:density}), we have
\begin{eqnarray*}
f_n\left(x\right)-f\left(x\right)&=&\left(1-\gamma_n\right)\left(f_{n-1}\left(x\right)-f\left(x\right)\right)+\gamma_n\left(Z_n\left(x\right)-f\left(x\right)\right)\\
&=&\sum_{k=1}^{n-1}\left[\prod_{j=k+1}^n\left(1-\gamma_j\right)\right]\gamma_k\left(Z_k\left(x\right)-f\left(x\right)\right)+\gamma_n\left(Z_n\left(x\right)-f\left(x\right)\right)+\left[\prod_{j=1}^n\left(1-\gamma_j\right)\right]\left(f_0\left(x\right)-f\left(x\right)\right)\\
&=&\Pi_n\sum_{k=1}^n\Pi_k^{-1}\gamma_k\left(Z_k\left(x\right)-f\left(x\right)\right)+\Pi_n\left(f_0\left(x\right)-f\left(x\right)\right).
\end{eqnarray*}
It follows that
\begin{eqnarray*}
\mathbb{E}\left[f_n\left(x\right)\right]-f\left(x\right)=\Pi_n\sum_{k=1}^n\Pi_k^{-1}\gamma_k\left(\mathbb{E}\left[Z_k\left(x\right)\right]-f\left(x\right)\right)+\Pi_n\left(f_0\left(x\right)-f\left(x\right)\right).
\end{eqnarray*}
Then, we can write that
\begin{eqnarray*}
f_n\left(x\right)-\mathbb{E}\left[f_n\left(x\right)\right]&=&\Pi_n\sum_{k=1}^n\Pi_k^{-1}\gamma_k\left(Z_k\left(x\right)-\mathbb{E}\left[Z_k\left(x\right)\right]\right)\\
&=&\Pi_n\sum_{k=1}^n\Pi_k^{-1}\gamma_kh_k^{-d}\left(Y_k-\mathbb{E}\left[Y_k\right]\right)
\end{eqnarray*}
Let $\left(\Psi_n\right)$ and $\left(B_n\right)$ be the sequences defined as
\begin{eqnarray*}
\Psi_n\left(x\right)&=&\Pi_n\sum_{k=1}^n\Pi_k^{-1}\gamma_kh_k^{-d}\left(Y_k-\mathbb{E}\left[Y_k\right]\right)\\
B_n\left(x\right)&=&\mathbb{E}\left[f_n\left(x\right)\right]-f\left(x\right)
\end{eqnarray*}
We have:
\begin{eqnarray}\label{eq:psiB}
f_n\left(x\right)-f\left(x\right)=\Psi_n\left(x\right)+B_n\left(x\right)
\end{eqnarray}
Theorems~\ref{the:LDPWW},~\ref{the:LDPDeu},~\ref{the:MDP},~\ref{theo:unifbounded} and~\ref{theo:unifunbounded} are consequences of~(\ref{eq:psiB}) and the following propositions.
\begin{prop}[{Pointwise LDP and MDP for $\left(\Psi_n\right)$}]\label{prop:LMDP}\quad \quad \quad \quad \quad \quad \quad \quad \quad \quad \quad \quad \quad \quad \quad \quad \quad \quad \quad \quad
\begin{enumerate}
\item Under the assumptions $\left(L1\right)$ and $\left(L2\right)$, the sequence $\left(f_n\left(x\right)-\mathbb{E}\left(f_n\left(x\right)\right)\right)$ satisfies a LDP with speed $\left(nh_n^d\right)$ and rate function $I_{a,x}$.
\item Under the assumptions $\left(L1\right)$ and $\left(L3\right)$, the sequence $\left(f_n\left(x\right)-\mathbb{E}\left(f_n\left(x\right)\right)\right)$ satisfies a LDP with speed $\left(nh_n^d\right)$ and rate function $I_{x}$.
\item Under the assumptions $\left(M1\right)-\left(M4\right)$, the sequence $\left(v_n\Psi_n\left(x\right)\right)$ satisfies a LDP with speed $\left(h_n^d/\left(\gamma_nv_n^2\right)\right)$ and rate function $J_{a,\alpha,x}$.
\end{enumerate}
\end{prop}
\begin{prop}[{Uniform LDP and MDP for $\left(\Psi_n\right)$}]\label{pr:psiunif}
\quad \quad \quad \quad \quad \quad \quad \quad \quad \quad \quad \quad \quad \quad \quad \quad \quad \quad \quad \quad
\begin{enumerate}
\item Let $\left(U1\right)-\left(U3\right)$ hold. Then for any bounded subset $U$ of $\mathbb{R}^d$ and for all $\delta>0$,
\begin{eqnarray*}
\lim_{n\to \infty} \frac{\gamma_nv_n^2}{h_n^d}\log \mathbb{P}\left[\sup_{x\in U}v_n\left|\Psi_n\left(x\right)\right|\geq \delta\right]=-\tilde{g}_U\left(\delta\right)
\end{eqnarray*} 
\item Let $\left(U1\right)-\left(U6\right)$ hold. Then for any subset $U$ of $\mathbb{R}^d$ and for all $\delta>0$,
\begin{eqnarray*}
-\tilde{g}_U\left(\delta\right)&\leq &\liminf_{n\to \infty} \frac{\gamma_nv_n^2}{h_n^d}\log \mathbb{P}\left[\sup_{x\in U}v_n\left|\Psi_n\left(x\right)\right|\geq \delta\right]\\
&\leq&\limsup_{n\to \infty} \frac{\gamma_nv_n^2}{h_n^d}\log \mathbb{P}\left[\sup_{x\in U}v_n\left|\Psi_n\left(x\right)\right|\geq \delta\right]\leq -\frac{\xi}{\xi+d}\tilde{g}_U\left(\delta\right)
\end{eqnarray*} 
\end{enumerate}
\end{prop}
The proof of the following proposition is given in Mokkadem et al. (2009).
\begin{prop}[{Pointwise and uniform convergence rate of $\left(B_n\right)$}]\label{pr:Bn}
$ $\\
Let Assumptions $\left( M1\right)-\left( M3\right)$ hold.
\begin{enumerate}
\item If for all $i,j\in\{1,\ldots d\}$, $\partial^2f/\partial x_i\partial x_j$ is continuous at $x$. We have
\begin{description}
\item If $a\leq\alpha/(d+4)$, then
\begin{eqnarray*}
B_n\left(x\right)=O\left(h_n^2\right).
\end{eqnarray*}
\item If $a>\alpha/(d+4)$, then
\begin{eqnarray*}
B_n\left(x\right)=o\left(\sqrt{\gamma_nh_n^{-d}}\right).
\end{eqnarray*}
\end{description}
\item If $\left(U2\right)$ holds, then:
\begin{description}
\item If $a\leq\alpha/(d+4)$, then
\begin{eqnarray*}
\sup_{x\in \mathbb{R}^d}\left|B_n\left(x\right)\right|=O\left(h_n^2\right).
\end{eqnarray*}
\item If $a>\alpha/(d+4)$, then
\begin{eqnarray*}
\sup_{x\in \mathbb{R}^d}\left|B_n\left(x\right)\right|=o\left(\sqrt{\gamma_nh_n^{-d}}\right).
\end{eqnarray*}
\end{description}
\end{enumerate}
\end{prop}
Set $x\in \mathbb{R}^d$; since the assumptions of Theorems~\ref{the:LDPWW} and~\ref{the:LDPDeu} guarantee that $\lim_{n\to \infty}B_n\left(x\right)=0$, Theorem~\ref{the:LDPWW} (respectively Theorem~\ref{the:LDPDeu}) is a straightforward consequence of the application of Part 1 (respectively of Part 2) of Proposition~\ref{prop:LMDP}. Moreover, under the assumptions of Theorem~\ref{the:MDP}, we have by application of Propostion~\ref{pr:Bn}, $\lim_{n\to \infty}v_nB_n\left(x\right)=0$; Theorem~\ref{the:MDP} thus straightfully follows from the application of Part 3 of Proposition~\ref{prop:LMDP}. Finaly, Theorem~\ref{theo:unifbounded} and~\ref{theo:unifunbounded} follows from Proposition~\ref{pr:psiunif} and the second part of Proposition~\ref{pr:Bn}.\\
We now state a preliminary lemma, which will be used in the proof of Proposition~\ref{prop:LMDP}.\\
For any $u\in \mathbb{R}$, Set
\begin{eqnarray*}
\Lambda_{n,x}\left(u\right)&=&\frac{\gamma_nv_n^2}{h_n^d}\log \mathbb{E}\left[\exp\left(u\frac{h_n^d}{\gamma_nv_n}\Psi_n\left(x\right)\right)\right]\\
\Lambda_{x}^{L,1}\left(u\right)&=&f\left(x\right)\left(\psi_a\left(u\right)-u\right),\\
\Lambda_{x}^{L,2}\left(u\right)&=&f\left(x\right)\left(\psi\left(u\right)-u\right),\\
\Lambda_{x}^M\left(u\right)&=&\frac{u^2}{2\left(2-\left(\alpha-ad\right)\xi\right)}f\left(x\right)
\int_{\mathbb{R}^d}K^2\left(z\right)dz
\end{eqnarray*}
\begin{lemma}\label{lemma:convLam}[Convergence of $\Lambda_{n,x}$]
\begin{enumerate}
\item (Pointwise convergence)\\
 If $f$ is continuous at $x$, then for all $u\in \mathbb{R}$
\begin{eqnarray}\label{con:lamd}
\lim_{n\to \infty}\Lambda_{n,x}\left(u\right)=\Lambda_{x}\left(u\right)
\end{eqnarray}
where
\begin{eqnarray*}
\Lambda_x\left(u\right)=\left\{\begin{array}{lllll}
\Lambda_x^{L,1}\left(u\right) & \mbox{when} & v_n\equiv 1& \mbox{,}&\left(L1\right)\,\,\mbox{and}\,\,\left(L2\right)\,\,\mbox{hold}\\
\Lambda_x^{L,2}\left(u\right) & \mbox{when} & v_n\equiv 1&\mbox{,}&
\left(L1\right)\,\,\mbox{and}\,\,\left(L3\right)\,\,\mbox{hold}\\
\Lambda_x^{M}\left(u\right) & \mbox{when} & v_n\to \infty & \mbox{,}&\left(M1\right)-\left(M4\right)\,\,\mbox{hold} \\
\end{array}
\right.
\end{eqnarray*}
\item (Uniform convergence)\\
If $f$ is uniformly continuous, then the convergence~(\ref{con:lamd}) holds uniformly in $x\in U$.
\end{enumerate}
\end{lemma}
Our proofs are now organized as follows: Lemma~\ref{lemma:convLam} is proved in Section~\ref{proof:lemma1}, Proposition~\ref{prop:LMDP} in Section~\ref{proof:prop2} and Proposition~\ref{pr:psiunif} in Section~\ref{proof:prop3}.
\subsection{Proof of Lemma~\ref{lemma:convLam}.}\label{proof:lemma1}
Set $u\in \mathbb{R}$, $u_n=u/v_n$ and $a_n=h_n^d\gamma_n^{-1}$. We have:
\begin{eqnarray*}
\Lambda_{n,x}\left(u\right)&=&\frac{v_n^2}{a_n}\log \mathbb{E}\left[\exp\left(u_na_n\Psi_n\left(x\right)\right)\right]\\
&=&\frac{v_n^2}{a_n}\log \mathbb{E}\left[\exp\left(u_na_n\Pi_n\sum_{k=1}^n\Pi_k^{-1}a_k^{-1}\left(Y_k-\mathbb{E}\left[Y_k\right]\right)\right)\right]\\
&=&\frac{v_n^2}{a_n}\sum_{k=1}^n\log \mathbb{E}\left[\exp\left(u_n\frac{a_n\Pi_n}{a_k\Pi_k}Y_k\right)\right]
-uv_n \Pi_n\sum_{k=1}^n\Pi_k^{-1}a_k^{-1}\mathbb{E}\left[Y_k\right]
\end{eqnarray*}
By Taylor expansion, there exists $c_{k,n}$ between $1$ and $\mathbb{E}\left[\exp\left(u_n\frac{a_n\Pi_n}{a_k\Pi_k}Y_k\right)\right]$ such that 
\begin{eqnarray*}
\log \mathbb{E}\left[\exp\left(u_n\frac{a_n\Pi_n}{a_k\Pi_k}Y_k\right)\right]=\mathbb{E}\left[\exp\left(u_n\frac{a_n\Pi_n}{a_k\Pi_k}Y_k\right)-1\right]
-\frac{1}{2c_{k,n}^2}\left(\mathbb{E}\left[\exp\left(u_n\frac{a_n\Pi_n}{a_k\Pi_k}Y_k\right)-1\right]\right)^2
\end{eqnarray*}
and $\Lambda_{n,x}$ can be rewriten as
\begin{eqnarray}\label{eq:lamb}
\Lambda_{n,x}\left(u\right)&=&\frac{v_n^2}{a_n}\sum_{k=1}^n\mathbb{E}\left[\exp\left(u_n\frac{a_n\Pi_n}{a_k\Pi_k}Y_k\right)-1\right]
-\frac{v_n^2}{2a_n}\sum_{k=1}^n
\frac{1}{c_{k,n}^2}\left(\mathbb{E}\left[\exp\left(u_n\frac{a_n\Pi_n}{a_k\Pi_k}Y_k\right)-1\right]\right)^2\nonumber\\
&&-uv_n \Pi_n\sum_{k=1}^n\Pi_k^{-1}a_k^{-1}\mathbb{E}\left[Y_k\right]
\end{eqnarray}
\paragraph{First case: $v_n\to \infty$.} A Taylor's expansion implies the existence of $c^{\prime}_{k,n}$ between $0$ and $u_n\frac{a_n\Pi_n}{a_k\Pi_k}Y_k$ such that 
\begin{eqnarray*}
\mathbb{E}\left[\exp\left(u_n\frac{a_n\Pi_n}{a_k\Pi_k}Y_k\right)-1\right]=u_n\frac{a_n\Pi_n}{a_k\Pi_k}\mathbb{E}\left[Y_k\right]+\frac{1}{2}\left(u_n\frac{a_n\Pi_n}{a_k\Pi_k}\right)^2\mathbb{E}\left[Y_k^2\right]+\frac{1}{6}\left(u_n\frac{a_n\Pi_n}{a_k\Pi_k}\right)^3\mathbb{E}\left[Y_k^3e^{c^{\prime}_{k,n}}\right]
\end{eqnarray*}
Therefore,
\begin{eqnarray}\label{eq:Moderate}
\Lambda_{n,x}\left(u\right)&=&\frac{1}{2}u^2a_n\Pi_n^2\sum_{k=1}^n\Pi_k^{-2}a_k^{-2}\mathbb{E}\left[Y_k^2\right]+\frac{1}{6}u^2u_na_n^2\Pi_n^3\sum_{k=1}^n\Pi_k^{-3}a_k^{-3}\mathbb{E}\left[Y_k^3e^{c^{\prime}_{k,n}}\right]\nonumber\\
&&-\frac{v_n^2}{2a_n}\sum_{k=1}^n
\frac{1}{c_{k,n}^2}\left(\mathbb{E}\left[\exp\left(u_n\frac{a_n\Pi_n}{a_k\Pi_k}Y_k\right)-1\right]\right)^2\nonumber\\
&=&\frac{1}{2}f\left(x\right)u^2a_n\Pi_n^2\sum_{k=1}^n\Pi_k^{-2}a_k^{-1}\gamma_k\int_{\mathbb{R}^d}K^2\left(z\right)dz+R_{n,x}^{\left(1\right)}\left(u\right)+R_{n,x}^{\left(2\right)}\left(u\right)
\end{eqnarray}
with
\begin{eqnarray*}
R_{n,x}^{\left(1\right)}\left(u\right)&=&\frac{1}{2}u^2a_n\Pi_n^2\sum_{k=1}^n\Pi_k^{-2}a_k^{-1}\gamma_k\int_{\mathbb{R}^d}K^2\left(z\right)\left[f\left(x-zh_k\right)-f\left(x\right)\right]dz\\
R_{n,x}^{\left(2\right)}\left(u\right)&=&\frac{1}{6}\frac{u^3}{v_n}a_n^2\Pi_n^3\sum_{k=1}^n\Pi_k^{-3}a_k^{-3}\mathbb{E}\left[Y_k^3e^{c^{\prime}_{k,n}}\right]-\frac{v_n^2}{2a_n}\sum_{k=1}^n
\frac{1}{c_{k,n}^2}\left(\mathbb{E}\left[\exp\left(u_n\frac{a_n\Pi_n}{a_k\Pi_k}Y_k\right)-1\right]\right)^2
\end{eqnarray*}
Since $f$ is continuous, we have $\lim_{k\to \infty}\left|f\left(x-zh_k\right)-f\left(x\right)\right|=0$, and thus, by the dominated convergence theorem, $\left(M1\right)$ implies that
\begin{eqnarray*}
\lim_{k\to \infty} \int_{\mathbb{R}^d}K^2\left(z\right)\left|f\left(x-zh_k\right)-f\left(x\right)\right|dz=0.
\end{eqnarray*}
Since $\left(a_n\right)\in\mathcal{GS}\left(\alpha-ad\right)$, and $\lim_{n\to \infty}\left(n\gamma_n\right)>\left(\alpha-ad\right)/2$. Lemma~\ref{lemma:Tech} then ensures that
\begin{eqnarray}\label{eq:a2}
a_n\Pi_n^2\sum_{k=1}^n\Pi_k^{-2}a_k^{-1}\gamma_k=\frac{1}{\left(2-\left(\alpha-ad\right)\xi\right)}+o\left(1\right),
\end{eqnarray}
it follows that $\lim_{n\to \infty}\left|R_{n,x}^{\left(1\right)}\left(u\right)\right|=0$.\\
Moreover, in view of~(\ref{eq:Y}), we have $\left|Y_k\right|\leq \left\|K\right\|_{\infty}$, then
\begin{eqnarray}\label{Yk:K}
c^{\prime}_{k,n}&\leq &\left|u_n\frac{a_n\Pi_n}{a_k\Pi_k}Y_k\right|\nonumber\\
&\leq & \left|u_n\right|\left\|K\right\|_{\infty}
\end{eqnarray}
Noting that $\mathbb{E}\left|Y_k\right|^3\leq h_k^d\left\|f\right\|_{\infty} \int_{\mathbb{R}^d}\left|K^3\left(z\right)\right|dz$. Hence, using Lemma~\ref{lemma:Tech} and~(\ref{Yk:K}), there exists a positive constant $c_1$ such that, for $n$ large enough,

\begin{eqnarray}\label{eq:R2term1}
\left|\frac{u^3}{v_n}a_n^2\Pi_n^3\sum_{k=1}^n\Pi_k^{-3}a_k^{-3}\mathbb{E}\left[Y_k^3e^{c^{\prime}_{k,n}}\right]\right|&\leq& c_1e^{\left|u_n\right|\left\|K\right\|_{\infty}}\frac{u^3}{v_n}\left\|f\right\|_{\infty}\int_{\mathbb{R}^d}\left|K^3\left(z\right)\right|dz
\end{eqnarray}
which goes to $0$ as $n\to \infty$ since $v_n\to \infty$.\\
Moreover, Lemma~\ref{lemma:Tech} ensures that
\begin{eqnarray}\label{eq:R2term2}
\lefteqn{\left|\frac{v_n^2}{2a_n}\sum_{k=1}^n
\frac{1}{c_{k,n}^2}\left(\mathbb{E}\left[\exp\left(u_n\frac{a_n\Pi_n}{a_k\Pi_k}Y_k\right)-1\right]\right)^2\right|}\nonumber\\
&& \leq  \frac{v_n^2}{2a_n}\sum_{k=1}^n
\left(\mathbb{E}\left[\exp\left(u_n\frac{a_n\Pi_n}{a_k\Pi_k}Y_k\right)-1\right]\right)^2\nonumber\\
&& \leq  \frac{u^2}{2}\left\|f\right\|_{\infty}^2a_n\Pi_n^2\sum_{k=1}^n
\Pi_k^{-2}a_k^{-1}\gamma_kh_k^d +o\left(a_n\Pi_n^2\sum_{k=1}^n
\Pi_k^{-2}a_k^{-1}\gamma_kh_k^d\right)\nonumber\\
&&=o\left(1\right)
\end{eqnarray}
The combination of~(\ref{eq:R2term1}) and~(\ref{eq:R2term2}) ensures that   
 $\lim_{n\to \infty}\left|R_{n,x}^{\left(2\right)}\left(u\right)\right|=0$. Then, we obtain from (\ref{eq:Moderate}) and (\ref{eq:a2}), $\lim_{n\to \infty}\Lambda_{n,x}\left(u\right)=\Lambda_x^{M}\left(u\right)$. 
\paragraph{Second case: $\left(v_n\right)\equiv 1$.}
It follows from~(\ref{eq:lamb}) that
\begin{eqnarray}\label{eq:lamb1}
\Lambda_{n,x}\left(u\right)&=&\frac{1}{a_n}\sum_{k=1}^n\mathbb{E}\left[\exp\left(u\frac{a_n\Pi_n}{a_k\Pi_k}Y_k\right)-1\right]
-\frac{1}{2a_n}\sum_{k=1}^n
\frac{1}{c_{k,n}^2}\left(\mathbb{E}\left[\exp\left(u\frac{a_n\Pi_n}{a_k\Pi_k}Y_k\right)-1\right]\right)^2\nonumber\\
&&-u \Pi_n\sum_{k=1}^n\Pi_k^{-1}a_k^{-1}\mathbb{E}\left[Y_k\right]\nonumber\\
&=&\frac{1}{a_n}\sum_{k=1}^nh_k^d\int_{\mathbb{R}^d}\left[\exp\left(u\frac{a_n\Pi_n}{a_k\Pi_k}K\left(z\right)\right)-1\right]f\left(x\right)dz-u \Pi_n\sum_{k=1}^n\Pi_k^{-1}\gamma_k\int_{\mathbb{R}^d}K\left(z\right)f\left(x\right)dz\nonumber\\
&&-R_{n,x}^{\left(3\right)}\left(u\right)+R_{n,x}^{\left(4\right)}\left(u\right)\nonumber\\
&=&f\left(x\right)\frac{1}{a_n}\sum_{k=1}^nh_k^d\left[\int_{\mathbb{R}^d}\left(\exp\left(uV_{n,k}K\left(z\right)\right)-1\right)-u V_{n,k}K\left(z\right)\right]dz\nonumber\\
&&-R_{n,x}^{\left(3\right)}\left(u\right)+R_{n,x}^{\left(4\right)}\left(u\right)
\end{eqnarray}
with 
\begin{eqnarray*}
V_{n,k}&=&\frac{a_n\Pi_n}{a_k\Pi_k}\\
R_{n,x}^{\left(3\right)}\left(u\right)&=&\frac{1}{2a_n}\sum_{k=1}^n
\frac{1}{c_{k,n}^2}\left(\mathbb{E}\left[\exp\left(u\frac{a_n\Pi_n}{a_k\Pi_k}Y_k\right)-1\right]\right)^2\nonumber\\
R_{n,x}^{\left(4\right)}\left(u\right)&=&\frac{1}{a_n}\sum_{k=1}^nh_k^d\int_{\mathbb{R}^d}\left[\exp\left(u\frac{a_n\Pi_n}{a_k\Pi_k}K\left(z\right)\right)-1\right]\left[f\left(x-zh_k\right)-f\left(x\right)\right]dz\\
&&-u \Pi_n\sum_{k=1}^n\Pi_k^{-1}\gamma_k\int_{\mathbb{R}^d}K\left(z\right)\left[f\left(x-zh_k\right)-f\left(x\right)\right]dz.
\end{eqnarray*}
 It follows from~(\ref{eq:R2term2}), that $\lim_{n\to \infty}\left|R^{\left(3\right)}_{n,x}\left(u\right)\right|=0$. \\
Since $\left|e^{t}-1\right|\leq \left|t\right|e^{\left|t\right|}$, we have 
\begin{eqnarray*}
\left|R_{n,x}^{\left(4\right)}\left(u\right)\right|
&\leq&\frac{1}{a_n}\sum_{k=1}^nh_k^d\int_{\mathbb{R}^d}\left|\left[\exp\left(u\frac{a_n\Pi_n}{a_k\Pi_k}K\left(z\right)\right)-1\right]\left[f\left(x-zh_k\right)-f\left(x\right)\right]\right|dz\\
&&+\left|u\right| \Pi_n\sum_{k=1}^n\Pi_k^{-1}\gamma_k\int_{\mathbb{R}^d}\left|K\left(z\right)\right|\left|f\left(x-zh_k\right)-f\left(x\right)\right|dz\\
&\leq&\left|u\right|e^{\left|u\right|\left\|K\right\|_{\infty}}\Pi_n\sum_{k=1}^n\Pi_k^{-1}\gamma_k\int_{\mathbb{R}^d}\left|K\left(z\right)\right|\left|f\left(x-zh_k\right)-f\left(x\right)\right|dz\\
&&+\left|u\right| \Pi_n\sum_{k=1}^n\Pi_k^{-1}\gamma_k\int_{\mathbb{R}^d}\left|K\left(z\right)\right|\left|f\left(x-zh_k\right)-f\left(x\right)\right|dz\\
&\leq&\left|u\right|\left(e^{\left|u\right|\left\|K\right\|_{\infty}}+1\right)\Pi_n\sum_{k=1}^n\Pi_k^{-1}\gamma_k\int_{\mathbb{R}^d}\left|K\left(z\right)\right|\left|f\left(x-zh_k\right)-f\left(x\right)\right|dz
\end{eqnarray*}
In view of Lemma~\ref{lemma:Tech} the sequence $\left(\Pi_n\sum_{k=1}^n\Pi_k^{-1}\gamma_k\right)$ is bounded, then, the dominated convergence theorem ensures that $\lim_{n\to \infty}R_{n,x}^{\left(4\right)}\left(u\right)=0$.\\
In the case $f$ is uniformly continuous, set $\varepsilon>0$ and let $M>0$ such that $2\left\|f\right\|_{\infty}\int_{\left\|z\right\|\leq M}\left|K\left(z\right)\right|dz\leq \varepsilon/2$. We need to prove that for $n$ sufficiently large
$$
\sup_{x\in \mathbb{R}^d}\int_{\left\|z\right\|\leq M}\left|K\left(z\right)\right|\left|f\left(x-zh_k\right)-f\left(x\right)\right|dz\leq \varepsilon/2$$
which is a straightforward consequence of the uniform continuity of $f$.\\

Then, it follows from~(\ref{eq:lamb1}), that
\begin{eqnarray}\label{eq:lambdavnk}
\lim_{n\to \infty}\Lambda_{n,x}\left(u\right)&=&\lim_{n\to \infty}f\left(x\right)\frac{\gamma_n}{h_n^d}\sum_{k=1}^nh_k^d\int_{\mathbb{R}^d}\left[\left(\exp\left(uV_{n,k}K\left(z\right)\right)-1\right)-u V_{n,k}K\left(z\right)\right]dz
\end{eqnarray}
\paragraph{In the case when $\left(v_n\right)\equiv 1$,  $\left(L1\right)$ and $\left(L2\right)$ hold}
 
 $ $\\
We have
\begin{eqnarray*}
\frac{\Pi_n}{\Pi_k}&=&\prod_{j=k+1}^n\left(1-\gamma_j\right)\\
&=&\frac{k}{n},
\end{eqnarray*}
then, 
\begin{eqnarray*}
V_{n,k}&=&\frac{a_n\Pi_n}{a_k\Pi_k}\\
&=&\left(\frac{k}{n}\right)^{ad}.
\end{eqnarray*}
Consequently, it follows from~(\ref{eq:lambdavnk}) and from some analysis considerations that
\begin{eqnarray*}
\lim_{n\to \infty}\Lambda_{n,x}\left(u\right)&=&f\left(x\right)\int_{\mathbb{R}^d} \left[\int_{0}^1s^{-ad}\left(\exp\left(us^{ad}K\left(z\right)\right)-1-u s^{ad}K\left(z\right)\right)ds\right]dz\\
&=&\Lambda_{x}^{L,1}\left(u\right)
\end{eqnarray*}
\paragraph{In the case when $\left(v_n\right)\equiv 1$,  $\left(L1\right)$ and $\left(L3\right)$ hold}  $ $ \\
We have
\begin{eqnarray*}
\frac{\Pi_n}{\Pi_k}&=&\prod_{j=k+1}^n\left(1-\gamma_j\right)\\
&=&\prod_{j=k+1}^n\left(1-\frac{h_j^d}{\sum_{l=1}^jh_l^d}\right)\\
&=&\prod_{j=k+1}^n\frac{\sum_{l=1}^{j-1}h_l^d}{\sum_{l=1}^jh_l^d}\\
&=&\frac{\sum_{l=1}^{k}h_l^d}{\sum_{l=1}^nh_l^d}\\
&=&\frac{\sum_{l=1}^{k}h_l^d}{h_k^d}\frac{h_k^d}{h_n^d}\frac{h_n^d}{\sum_{l=1}^nh_l^d}\\
&=&\frac{\gamma_n}{\gamma_k}\frac{h_k^d}{h_n^d},
\end{eqnarray*}
then,
\begin{eqnarray*}
V_{n,k}=1.
\end{eqnarray*}
Consequently, it follows from~(\ref{eq:lambdavnk}) that
\begin{eqnarray*}
\lim_{n\to \infty}\Lambda_{n,x}\left(u\right)&=&
f\left(x\right)\int_{\mathbb{R}^d}
\left[\left(\exp\left(uK\left(z\right)\right)-1\right)-u K\left(z\right)\right]dz\\
&=&\Lambda_{x}^{L,2}\left(u\right)
\end{eqnarray*}
and thus Lemma 1 is proved.
\subsection{Proof of Proposition~\ref{prop:LMDP}}\label{proof:prop2}
$ $\\
To prove Proposition~\ref{prop:LMDP}, we apply Proposition 1 in Mokkadem et al. (2006), Lemma~\ref{lemma:convLam} and the following result (see Puhalskii, 1994).
\begin{lemma}\label{lemma:puhalskii}
Let $\left(Z_n\right)$ be a sequence of real random variables, $\left(\nu_n\right)$ a positive sequence satisfying $\lim_{n\to \infty}\nu_n=+\infty$, and suppose that there exists some convex non-negative function $\Gamma$ defined on $\mathbb{R}$ such that
\begin{eqnarray*}
\forall u\in \mathbb{R}, \lim_{n\to \infty}\frac{1}{\nu_n}\log \mathbb{E}\left[\exp\left(u\nu_nZ_n\right)\right]=\Gamma\left(u\right).
\end{eqnarray*}
If the Legendre function $\Gamma^*$ of $\Gamma$ is a strictly convex function, then the sequence $\left(Z_n\right)$ satisfies a \texttt{LDP} of speed $\left(\nu_n\right)$ and good rate fonction $\Gamma^*$.
\end{lemma}
In our framework, when $v_n\equiv 1$ and $\gamma_n=n^{-1}$, we take $Z_n=f_n\left(x\right)-\mathbb{E}\left(f_n\left(x\right)\right)$, $\nu_n=nh_n^d$ with $h_n=cn^{-a}$ where $a\in \left]0,1/d\right[$ and $\Gamma=\Lambda_x^{L,1}$. In this case, the Legendre transform of $\Gamma=\Lambda_x^{L,1}$ is the rate function $I_{a,x}:t\to f\left(x\right)I_a\left(\frac{t}{f\left(x\right)}+1\right)$ which is strictly convex by Proposition 1 in Mokkadem et al. (2006). Farther, when $v_n\equiv 1$ and $\gamma_n=h_n^d\left(\sum_{k=1}^nh_k^d\right)^{-1}$, we take $Z_n=f_n\left(x\right)-\mathbb{E}\left(f_n\left(x\right)\right)$, $\nu_n=nh_n^d$ with $h_n\in \mathcal{GS}\left(-a\right)$ where $a\in \left]0,1/d\right[$ and $\Gamma=\Lambda_x^{L,2}$. In this case, the Legendre transform of $\Gamma=\Lambda_x^{L,2}$ is the rate function $I_{x}:t\to f\left(x\right)I\left(\frac{t}{f\left(x\right)}+1\right)$ which is strictly convex by Proposition 1 in Mokkadem et al. (2005). Otherwise, when, $v_n\to \infty$, we take $Z_n=v_n\left(f_n\left(x\right)-\mathbb{E}\left(f_n\left(x\right)\right)\right)$, $\nu_n=h_n^d/\left(\gamma_nv_n^2\right)$ and $\Gamma=\Lambda_x^M$; $\Gamma^*$ is then the quadratic rate function $J_{a,\alpha,x}$ defined in~(\ref{sys:MDP}) and thus Proposition~\ref{prop:LMDP} follows.
\subsection{Proof of Proposition~\ref{pr:psiunif}}\label{proof:prop3}
In order to prove Proposition~\ref{pr:psiunif}, we first establish some lemmas.
\begin{lemma}\label{lem3}
Let $\phi:\mathbb{R}^+\to \mathbb{R}$ be the function defined for $\delta>0$ as
\begin{eqnarray*}
\phi\left(\delta\right)&=&\left\{\begin{array}{lllll}
\left(\psi_a^{\prime}\right)^{-1}\left(1+\frac{\delta}{\|f\|_{U,\infty}}\right) & \mbox{when} & v_n\equiv 1& \mbox{,}&\left(L1\right)\,\,\mbox{and}\,\,\left(L2\right)\,\,\mbox{hold}\\
\left(\psi^{\prime}\right)^{-1}\left(1+\frac{\delta}{\|f\|_{U,\infty}}\right) & \mbox{when} & v_n\equiv 1&\mbox{,}&\left(L1\right)\,\,\mbox{and}\,\,\left(L3\right)\,\,\mbox{hold}\\
\frac{\delta\left(2-\left(\alpha-ad\right)\xi\right)}{\|f\|_{U,\infty}\int_{\mathbb{R}^d}K^2\left(z\right)dz}
 & \mbox{when} & v_n\to \infty&\mbox{,}& \left(M1\right)-\left(M4\right)\,\,\mbox{hold} \\
\end{array}
\right.
\end{eqnarray*}
\begin{enumerate}
\item $\sup_{u\in \mathbb{R}}\left\{u\delta-\sup_{x\in U}\Lambda_x\left(u\right)\right\}$ equals $g_U\left(\delta\right)$ and is achieved for $u=\phi\left(\delta\right)>0$.
\item $\sup_{u\in \mathbb{R}}\left\{-u\delta-\sup_{x\in U}\Lambda_x\left(u\right)\right\}$ equals $g_U\left(\delta\right)$ and is achieved for $u=\phi\left(-\delta\right)<0$.
\end{enumerate}
\end{lemma}
\paragraph{Proof of Lemma~\ref{lem3}}.
We just prove the first part, the proof of the second part one being similar.
\begin{itemize}
\item First case : $v_n\equiv 1$, $\left(L1\right)$ and $\left(L2\right)$ hold.\\
Since $e^t\geq 1+t$, for all $t$, we have $\psi_a\left(u\right)\geq u$ and therefore,
\begin{eqnarray*}
u\delta-\sup_{x\in U}\Lambda_x\left(u\right)&=&u\delta-\|f\|_{U,\infty}\left(\psi_a\left(u\right)-u\right)\\
&=&\|f\|_{U,\infty}\left[u\left(1+\frac{\delta}{\|f\|_{U,\infty}}\right)-\psi_a\left(u\right)\right]
\end{eqnarray*}
The function $u\mapsto u\delta-\sup_{x\in U}\Lambda_x\left(u\right)$ has second derivative $-\|f\|_{U,\infty}\psi_a^{\prime\prime}\left(u\right)<0$ and thus it has a unique maximum achieved for 
\begin{eqnarray*}
u_0=\left(\psi_a^{\prime}\right)^{-1}\left(1+\frac{\delta}{\|f\|_{U,\infty}}\right)
\end{eqnarray*}
Now, since $\psi^{\prime}_a$ is increasing and since $\psi^{\prime}_a\left(0\right)=1$, we deduce that $u_0>0$.
\item Second case : $v_n\equiv 1$, $\left(L1\right)$ and $\left(L3\right)$ hold.\\
Since $e^t\geq 1+t$, for all $t$, we have $\psi\left(u\right)\geq u$ and therefore,
\begin{eqnarray*}
u\delta-\sup_{x\in U}\Lambda_x\left(u\right)&=&u\delta-\|f\|_{U,\infty}\left(\psi\left(u\right)-u\right)\\
&=&\|f\|_{U,\infty}\left[u\left(1+\frac{\delta}{\|f\|_{U,\infty}}\right)-\psi\left(u\right)\right]
\end{eqnarray*}
The function $u\mapsto u\delta-\sup_{x\in U}\Lambda_x\left(u\right)$ has second derivative $-\|f\|_{U,\infty}\psi^{\prime\prime}\left(u\right)<0$ and thus it has a unique maximum achieved for 
\begin{eqnarray*}
u_0=\left(\psi^{\prime}\right)^{-1}\left(1+\frac{\delta}{\|f\|_{U,\infty}}\right)
\end{eqnarray*}
Now, since $\psi^{\prime}$ is increasing and since $\psi^{\prime}\left(0\right)=1$, we deduce that $u_0>0$.
\item Third case $v_n\to \infty$ and $\left(M2\right)$ holds.
In this case, we have
\begin{eqnarray*}
u\delta-\sup_{x\in U}\Lambda_x\left(u\right)&=&u\delta-\frac{u^2}{2\left(2-\left(\alpha-ad\right)\xi\right)}\|f\|_{U,\infty}\int_{\mathbb{R}^d}K^2\left(z\right)dz.
\end{eqnarray*}
In view of the assumption $\left(M2\right)$, we have $\xi^{-1}>\left(\alpha-ad\right)/2$, then the function $u\mapsto u\delta-\sup_{x\in U}\Lambda_x\left(u\right)$ has second derivative $-\frac{1}{\left(2-\left(\alpha-ad\right)\xi\right)}\|f\|_{U,\infty}\int_{\mathbb{R}^d}K^^2\left(z\right)dz<0$ and thus it has a unique maximum achieved for 
\begin{eqnarray*}
u_0=\frac{\delta\left(2-\left(\alpha-ad\right)\xi\right)}{\|f\|_{U,\infty}\int_{\mathbb{R}^d}K^2\left(z\right)dz}>0
\end{eqnarray*}
\end{itemize}

\begin{lemma}\label{lem4}
 $ $\\
\begin{itemize}
\item In the case when $\left(v_n\right)\equiv 1$ and $\left(\gamma_n\right)=\left(n^{-1}\right)$, let $\left(L1\right)$ and $\left(L2\right)$ hold;
\item In the case when $\left(v_n\right)\equiv 1$ and $\left(\gamma_n\right)=\left(h_n^{d}\left(\sum_{k=1}^nh_k^d\right)^{-1}\right)$, let $\left(L1\right)$ and $\left(L3\right)$ hold;
\item In the case when $v_n\to \infty$, let $\left(M1\right)-\left(M4\right)$ hold.\\
Then for any $\delta>0$,
\begin{eqnarray*}
\lim_{n\to \infty}\frac{\gamma_nv_n^2}{h_n^d}\log \sup_{x\in U}\mathbb{P}\left[v_n\Psi_n\left(x\right)\geq \delta\right]&=&-g_U\left(\delta\right)\\
\lim_{n\to \infty}\frac{\gamma_nv_n^2}{h_n^d}\log \sup_{x\in U}\mathbb{P}\left[v_n\Psi_n\left(x\right)\leq -\delta\right]&=&-g_U\left(-\delta\right)\\
\lim_{n\to \infty}\frac{\gamma_nv_n^2}{h_n^d}\log \sup_{x\in U}\mathbb{P}\left[v_n\left|\Psi_n\left(x\right)\right|\leq -\delta\right]&=&-\tilde{g}_U\left(-\delta\right)
\end{eqnarray*}
\end{itemize}
\end{lemma}
\paragraph{Proof of Lemma~\ref{lem4}.} 
The proof of Lemma~\ref{lem4} is similar to the proof of Lemma 4 in Mokkadem et al. (2006).
\begin{lemma}\label{lem5}
Let Assumptions $\left(U1\right)-\left(U3\right)$ hold and assume that either $\left(v_n\right)\equiv 1$ or $\left(U4\right)$ holds.
\begin{enumerate}
\item If $U$ is a bounded set, then for any $\delta>0$, we have
\begin{eqnarray*}
\lim_{n\to \infty}\frac{\gamma_nv_n^2}{h_n^d}\log \mathbb{P}\left[\sup_{x\in U}v_n\left|\Psi_n\left(x\right)\right|\right]\leq -\tilde{g}_U\left(\delta\right)
\end{eqnarray*}
\item If $U$ is an unbounded set, then, for any $b>0$ and $\delta>0$,
\begin{eqnarray*}
\limsup_{n\to \infty}\frac{\gamma_nv_n^2}{h_n^d}\log \mathbb{P}\left[\sup_{x\in U,\|x\|\leq w_n}v_n\left|\Psi_n\left(x\right)\right|\right]\leq db-\tilde{g}_U\left(\delta\right)
\end{eqnarray*}
where $w_n=\exp\left(\frac{bh_n^d}{\gamma_nv_n^2}\right)$.
\end{enumerate}
\end{lemma}
\paragraph{Proof of Lemma~\ref{lem5}.}
Set $\rho\in \left]0,\delta\right[$, let $\beta$ denote the H\"older order of $K$, and $\|K\|_H$ its corresponding H\"older norm. Set $w_n=\exp\left(\frac{bh_n^d}{\gamma_nv_n^2}\right)$ and
\begin{eqnarray*}
R_n=\left(\frac{\rho}{2\|K\|_{H}v_n\Pi_n\sum_{k=1}^n\Pi_k^{-1}\gamma_kh_k^{-\left(d+\beta\right)}}\right)^{\frac{1}{\beta}}
\end{eqnarray*} 
We begin with the proof of the second part of Lemma~\ref{lem5}. There exist $N^{\prime}\left(n\right)$ points of $\mathbb{R}^d$, $y_1^{\left(n\right)},y_2^{\left(n\right)},\ldots,y_{N^{\prime}\left(n\right)}^{\left(n\right)}$ such that the ball $\left\{x\in \mathbb{R}^d; \|x\|\leq w_n\right\}$ can covered by the $N^{\prime}\left(n\right)$ balls $B_i^{\left(n\right)}=\left\{x\in \mathbb{R}^d; \|x-y_i^{\left(n\right)}\|\leq R_n\right\}$ and such that $N^{\prime}\left(n\right)\leq 2\left(\frac{2w_n}{R_n}\right)^d$. Considering only the $N\left(n\right)$ balls that intersect $\left\{x\in U; \|x\|\leq w_n\right\}$, we can write
\begin{eqnarray*}
\left\{x\in U; \|x\|\leq w_n\right\}\subset \cup_{i=1}^{N\left(n\right)}B_i^{\left(n\right)}.
\end{eqnarray*}
For each $i\in \left\{1,\ldots,N\left(n\right)\right\}$, set $x_i^{\left(n\right)}\in B_i^{\left(n\right)}\cap U$. We then have:
\begin{eqnarray*}
\mathbb{P}\left[\sup_{x\in U,\|x\|\leq w_n}v_n\left|\Psi_n\left(x\right)\right|\geq \delta\right]&\leq & \sum_{i=1}^{N\left(n\right)}\mathbb{P}\left[\sup_{x\in B_i^{\left(n\right)}}v_n\left|\Psi_n\left(x\right)\right|\geq \delta\right]\\
&\leq & N\left(n\right) \max_{1\leq i \leq N\left(n\right)} \mathbb{P}\left[\sup_{x\in B_i^{\left(n\right)}}v_n\left|\Psi_n\left(x\right)\right|\geq \delta\right].
\end{eqnarray*}
Now, for any $i\in \left\{1,\ldots,N\left(n\right)\right\}$ and any $x\in B_i^{\left(n\right)}$,
\begin{eqnarray*}
v_n\left|\Psi_n\right|&\leq & v_n \left|\Psi_n\left(x_i^{\left(n\right)}\right)\right|\\
&&+v_n\Pi_n\sum_{k=1}^n\Pi_k^{-1}\gamma_kh_k^{-d}\left|K\left(\frac{x-X_k}{h_k}\right)-K\left(\frac{x_i^{\left(n\right)}-X_k}{h_k}\right)\right|\\
&&+v_n\Pi_n\sum_{k=1}^n\Pi_k^{-1}\gamma_kh_k^{-d}\mathbb{E}\left|K\left(\frac{x-X_k}{h_k}\right)-K\left(\frac{x_i^{\left(n\right)}-X_k}{h_k}\right)\right|\\
&\leq & v_n \left|\Psi_n\left(x_i^{\left(n\right)}\right)\right|+2v_n\|K\|_H\Pi_n\sum_{k=1}^n\Pi_k^{-1}\gamma_kh_k^{-d}\left(\frac{\|x-x_i^{\left(n\right)}\|}{h_k}\right)^{\beta}\\
&\leq & v_n \left|\Psi_n\left(x_i^{\left(n\right)}\right)\right|+2v_n\|K\|_H\Pi_n\sum_{k=1}^n\Pi_k^{-1}\gamma_kh_k^{-\left(d+\beta\right)}R_n^{\beta}\\
&\leq & v_n \left|\Psi_n\left(x_i^{\left(n\right)}\right)\right|+\rho
\end{eqnarray*}
Hence, we deduce that
\begin{eqnarray*}
\mathbb{P}\left[\sup_{x\in U,\|x\|\leq w_n}v_n\left|\Psi_n\left(x\right)\right|\geq \delta \right]
&\leq &N\left(n\right)\max_{1\leq i\leq N\left(n\right)}
\mathbb{P}\left[v_n\left|\Psi_n\left(x_i^{\left(n\right)}\right)\right|\geq \delta-\rho\right]\\
&\leq & N\left(n\right)\sup_{x\in U}\mathbb{P}\left[v_n\left|\Psi_n\left(x_i^{\left(n\right)}\right)\right|\geq \delta-\rho\right]
\end{eqnarray*}
Further, by definition of $N\left(n\right)$ and $w_n$, we have 
\begin{eqnarray*}
\log N\left(n\right)\leq \log N^{\prime}\left(n\right) \leq db \frac{h_n^{d}}{\gamma_nv_n^2}+\left(d+1\right)\log 2-d \log R_n
\end{eqnarray*}
and 
\begin{eqnarray*}
\frac{\gamma_nv_n^2}{h_n^d}\log R_n = \frac{\gamma_nv_n^2}{\beta h_n^{d}}\left[\log \rho-\log \left(2\|K\|_H\right)-\log v_n -\log \left(\Pi_n\sum_{k=1}^n\Pi_k^{-1}\gamma_kh_k^{-\left(d+\beta\right)}\right)\right].
\end{eqnarray*}
Moreover, we have $\left(h_n^{\left(d+\beta\right)}\right)\in \mathcal{GS}\left(-a\left(d+\beta\right)\right)$. Lemma~\ref{lemma:Tech} ensures that
\begin{eqnarray*}
\Pi_n\sum_{k=1}^n\Pi_k^{-1}\gamma_kh_k^{-\left(d+\beta\right)}=O\left(h_n^{-\left(d+\beta\right)}\right),
\end{eqnarray*}
then, in view of $\left(U3\right)$, we have
\begin{eqnarray}\label{eq:17}
\limsup_{n\to \infty}\frac{\gamma_nv_n^2}{h_n^d}\log N\left(n\right)\leq db
\end{eqnarray}
The application of Lemma~\ref{lem4} then yiels
\begin{eqnarray*}
\limsup_{n\to \infty}\frac{\gamma_nv_n^2}{h_n^d}\log \mathbb{P}\left[\sup_{x\in U,\|x\|\leq w_n}v_n\left|\Psi_n\left(x\right)\right|\geq \delta\right]&\leq& \limsup_{n\to \infty}\frac{\gamma_nv_n^2}{h_n^d}\log N\left(n\right)-\tilde{g}_U\left(\delta-\rho\right)\\
&\leq & db-\tilde{g}_U\left(\delta-\rho\right).
\end{eqnarray*}
Since the inequality holds for any $\rho\in \left]0,\delta\right[$, part 2 of Lemma~\ref{lem5} thus follows from the continuity of $\tilde{g}_U$.\\

Let us now consider part 1 of Lemma~\ref{lem5}. This part is proved by following the same steps as for part 2, except that the number $N\left(n\right)$ of balls covering $U$ is at most the integer part of $\left(\Delta/R_n\right)^d$, where $\Delta$ denotes the diameter of $\overline{U}$. Relation~(\ref{eq:17}) then becomes 
\begin{eqnarray*}
\limsup_{n\to \infty}\frac{\gamma_nv_n^2}{h_n^{d}}\log R_n\leq 0
\end{eqnarray*}
and Lemma~\ref{lem5} is proved.
\begin{lemma}\label{lem6}
Let $\left(U1\right)i)$, $\left(M2\right)$ and $\left(U6\right)i)$ hold. Assume that either $\left(v_n\right)\equiv 1$ or $\left(U3\right)$ and $\left(U6\right)ii)$ hold. Moreover assume that $f$ is continuous. For any $b>0$ if we set $w_n=\exp\left(\frac{bh_n^d}{\gamma_nv_n^2}\right)$ then, for any $\rho>0$, we have, for $n$ large enough,
\begin{eqnarray*}
\sup_{x\in U, \|x\|\geq w_n}v_n\Pi_n\sum_{k=1}^n\Pi_k^{-1}\gamma_kh_k^{-d}\left|\mathbb{E}\left[K\left(\frac{x-X_k}{h_k}\right)\right]\right|\leq \rho
\end{eqnarray*}
\end{lemma}
\paragraph{Proof of Lemma~\ref{lem6}.}
We have 
\begin{eqnarray}\label{eq:18}
v_n\Pi_n\sum_{k=1}^n\Pi_k^{-1}\gamma_kh_k^{-d}\mathbb{E}\left[K\left(\frac{x-X_k}{h_k}\right)\right]=v_n\Pi_n\sum_{k=1}^n\Pi_k^{-1}\gamma_k\int_{\mathbb{R}^d}K\left(z\right)f\left(x-zh_k\right)dz.
\end{eqnarray}
First, Lemma~\ref{lemma:Tech}, ensures that 
\begin{eqnarray}\label{eq:lempi}
\Pi_n\sum_{k=1}^n\Pi_k^{-1}\gamma_k=1+o\left(1\right).
\end{eqnarray}
Set $\rho>0$. In the case $\left(v_n\right)\equiv 1$, we set $M$ such that $\|f_{\infty}\|\int_{\|z\|>M}\left|K\left(z\right)\right|dz\leq \rho/2$; it follows from~(\ref{eq:lempi}) that
\begin{eqnarray*}
\lefteqn{v_n\Pi_n\sum_{k=1}^n\Pi_k^{-1}\gamma_kh_k^{-d}\left|\mathbb{E}\left[K\left(\frac{x-X_k}{h_k}\right)\right]\right|}\\
&&\leq \frac{\rho}{2}+f\left(x\right)\int_{\|z\|\leq M}\left|K\left(z\right)\right|dz\\
&&+\Pi_n\sum_{k=1}^n\Pi_k^{-1}\gamma_k\int_{\|z\|>M}\left|K\left(z\right)\right|\left|f\left(x-zh_k\right)-f\left(x\right)\right|dz.
\end{eqnarray*}
Lemma~\ref{lem6} then follows from the fact that $f$ fulfills $\left(U6\right) ii)$. As matter of fact, this conditions implies that $\lim_{\|x\|\to \infty, x\in \overline{U}}f\left(x\right)=0$ and that the third term in the right-hand-side of the previous inequality goes to $0$ as $n\to \infty$ (by the dominated convergence).\\
Let us now assume that $\lim_{n\to \infty}v_n=\infty$; relation~(\ref{eq:18}) can be rewritten as
\begin{eqnarray*}
v_n\Pi_n\sum_{k=1}^n\Pi_k^{-1}\gamma_kh_k^{-d}\mathbb{E}\left[K\left(\frac{x-X_k}{h_k}\right)\right]&=&v_n\Pi_n\sum_{k=1}^n\Pi_k^{-1}\gamma_k\int_{\|z\|\leq w_n/2}K\left(z\right)f\left(x-zh_k\right)dz\\
&&+v_n\Pi_n\sum_{k=1}^n\Pi_k^{-1}\gamma_k\int_{\|z\|\geq w_n/2}K\left(z\right)f\left(x-zh_k\right)dz.
\end{eqnarray*} 
First, since $\|x\|\geq w_n$ and $\|z\|\leq w_n/2$, we have
\begin{eqnarray*}
\|x-zh_k\|&\geq&w_n\left(1-h_i/2\right)\\
&\geq & w_n/2\quad \mbox{for $n$ large enough}.
\end{eqnarray*}
Moreover, in view of assumptions $\left(U3\right)$, for all $\xi>0$, 
\begin{eqnarray}\label{eq:vnwn}
\lim_{n\to \infty}\frac{v_n}{w_n^{\xi}}=\lim_{n\to \infty}\exp\left\{-\xi b\frac{h_n^d}{\gamma_nv_n^2}\left(1-\frac{v_n^2\log v_n}{\xi bh_n^d}\right)\right\}=  0.
\end{eqnarray}
Set $M_f=\sup_{x\in \mathbb{R}^d}\|x\|^{\eta}f\left(x\right)$. Assumption $\left(U6\right) ii)$ and equations~(\ref{eq:lempi}),~(\ref{eq:vnwn}) implie that, for $n$ sufficiently large,
\begin{eqnarray*}
\lefteqn{\sup_{\|x\|\geq w_n}v_n\Pi_n\sum_{k=1}^n\Pi_k^{-1}\gamma_k\int_{\|z\|\leq w_n/2}\left|K\left(z\right)f\left(x-zh_k\right)\right|dz}\\
&&\leq M_f\sup_{\|x\|\geq w_n}v_n\Pi_n\sum_{k=1}^n\Pi_k^{-1}\gamma_k\int_{\|z\|\leq w_n/2}\left|K\left(z\right)\right|\|x-zh_k\|^{-\eta}dz\\
&&\leq 2^{\eta}M_f\frac{v_n}{w_n^{\eta}}\int_{\mathbb{R}^d}\left|K\left(z\right)\right|dz\\
&&\leq \frac{\rho}{2}.
\end{eqnarray*}
Moreover, in view of $\left(U3\right)$, $\left(U6\right)i)$ and~(\ref{eq:lempi}),~(\ref{eq:vnwn}), for $n$ sufficiently large, 
\begin{eqnarray*}
\lefteqn{\sup_{\|x\|\geq w_n}v_n\Pi_n\sum_{k=1}^n\Pi_k^{-1}\gamma_k\int_{\|z\|> w_n/2}\left|K\left(z\right)f\left(x-zh_k\right)\right|dz}\\
&&\leq 2^{\zeta}M_f\frac{v_n}{w_n^{\zeta}}\int_{\|z\|>w_n/2}\|z\|^{\zeta}\left|K\left(z\right)\right|dz\\
&&\leq \frac{\rho}{2}.
\end{eqnarray*}
This concludes the proof of Lemma~\ref{lem6}.
Since $K$ is a bounded function that vanishes at infinity, we have $\lim_{\|x\|\to \infty}\left|\Psi_n\left(x\right)\right|=0$ for every $n\geq 1$. Moreover, since $K$ is assumed to be continuous, $\Psi_n$ is continuous, and this ensures the existence of a random variable $s_n$ such that
\begin{eqnarray*}
\left|\Psi_n\left(s_n\right)\right|=\sup_{x\in U}\left|\Psi_n\left(x\right)\right|.
\end{eqnarray*}
\begin{lemma}\label{lem7}
$ $\\
Let Assumptions $\left(U1\right)-\left(U3\right)$, $\left(U4\right) ii)$ and $\left(U5\right)$ hold. Suppose either $\left(v_n\right)\equiv 1$ or $\left(H6\right)$ hold. For any $b>0$, set $w_n=\exp\left(b\frac{h_n^d}{\gamma_nv_n^2}\right)$; for any $\delta>0$, we have
\begin{eqnarray}
\limsup_{n\to \infty}\frac{\gamma_nv_n^2}{h_n^d}\log \mathbb{P}\left[\|s_n\|\geq w_n\quad \mbox{and}\quad \left|\Psi_n\left(s_n\right)\right|\geq \delta\right]\quad \leq \quad -b\beta
\end{eqnarray}
\end{lemma}
\paragraph{Proof of Lemma~\ref{lem7}.}
We first note that $s_n\in \overline{U}$ and therefore
\begin{eqnarray*}
\lefteqn{\|s_n\|\geq w_n \quad \mbox{and} \quad v_n\left|\Psi_n\left(s_n\right)\right|\geq \delta}\\
&&\Rightarrow \|s_n\|\geq w_n \quad \mbox{and} \quad v_n\left|\Pi_n\sum_{k=1}^n\Pi_k^{-1}\gamma_kh_k^{-d}K\left(\frac{s_n-X_k}{h_k}\right)\right|\\
&&\quad \quad +v_n\mathbb{E}\left|\Pi_n\sum_{k=1}^n\Pi_k^{-1}\gamma_kh_k^{-d}K\left(\frac{s_n-X_k}{h_k}\right)\right|\geq \delta\\
&&\Rightarrow \|s_n\|\geq w_n \quad \mbox{and} \quad v_n\Pi_n\sum_{k=1}^n\Pi_k^{-1}\gamma_kh_k^{-d}\left|K\left(\frac{s_n-X_k}{h_k}\right)\right|\\
&&\quad \quad -\sup_{\|x\|\geq w_n,x\in \overline{U}}v_n\Pi_n\sum_{k=1}^n\Pi_k^{-1}\gamma_kh_k^{-d}\mathbb{E}\left|K\left(\frac{s_n-X_k}{h_k}\right)\right|\geq \delta.
\end{eqnarray*}
Set $\rho\in \left]0,\delta\right[$; the application of Lemma~\ref{lem6} ensures that, for $n$ large enough, 
\begin{eqnarray*}
\lefteqn{\|s_n\|\geq w_n \quad \mbox{and} \quad v_n\left|\Psi_n\left(s_n\right)\right|\geq \delta}\\
&&\Rightarrow \|s_n\|\geq w_n \quad \mbox{and} \quad v_n\left|\Pi_n\sum_{k=1}^n\Pi_k^{-1}\gamma_kh_k^{-d}K\left(\frac{s_n-X_k}{h_k}\right)\right|\geq\delta-\rho.
\end{eqnarray*}
Set $\kappa=\sup_{x\in \mathbb{R}^d}\|x\|^{\gamma}\left|K\left(x\right)\right|$ (see Assumption $\left(U5\right)$). We obtain, for $n$ sufficiently large,
\begin{eqnarray*}
\lefteqn{\|s_n\|\geq w_n \quad \mbox{and} \quad v_n\left|\Psi_n\left(s_n\right)\right|\geq \delta}\\
&&\Rightarrow \|s_n\|\geq w_n \quad \mbox{and} \quad \quad \exists k\in \left\{1,\ldots,n\right\} \quad \mbox{such that}\quad  \frac{v_n}{h_k^d}\left|K\left(\frac{s_n-X_k}{h_k}\right)\right|\geq\delta-\rho\\
&&\Rightarrow \|s_n\|\geq w_n \quad \mbox{and} \quad \quad \exists k\in \left\{1,\ldots,n\right\} \quad \mbox{such that}\quad  \kappa h_k^{\gamma}\geq \frac{h_k^d}{v_n}\|s_n-X_k\|^{\gamma}\left(\delta-\rho\right)\\
&&\Rightarrow \|s_n\|\geq w_n \quad \mbox{and} \quad \quad \exists k\in \left\{1,\ldots,n\right\} \quad \mbox{such that}\quad  \left|\|s_n\|-\|X_k\|\right|\leq \left[\frac{\kappa v_nh_k^{\gamma-d}}{\delta-\rho}\right]^{\frac{1}{\gamma}}\\
&&\Rightarrow \|s_n\|\geq w_n \quad \mbox{and} \quad \quad \exists k\in \left\{1,\ldots,n\right\} \quad \mbox{such that}\quad  \|X_k\|\leq \|s_n\| -\left[\frac{\kappa v_nh_k^{\gamma-d}}{\delta-\rho}\right]^{\frac{1}{\gamma}}\\
&&\Rightarrow \|s_n\|\geq w_n \quad \mbox{and} \quad \quad \exists k\in \left\{1,\ldots,n\right\} \quad \mbox{such that}\quad  \|X_k\|\leq w_n\left(1-u_{n,k}\right)\quad \mbox{with}\\
&&\quad u_{n,k}=w_n^{-1}v_n^{\frac{1}{\gamma}}h_k^{\frac{\gamma-d}{\gamma}}\left(\frac{\kappa}{\delta-\rho}\right)^{\frac{1}{\gamma}}.
\end{eqnarray*}
Moreover, we can write $u_{n,k}$ as 
\begin{eqnarray*}
u_{n,k}=\exp\left(-b\frac{h_n^d}{\gamma_nv_n^2}\left[1-\frac{1}{b\gamma}\frac{\gamma_nv_n^2\log v_n}{h_n^{d}}-\frac{\gamma-d}{b\gamma}\frac{\gamma_nv_n^2\log\left(h_k\right)}{h_n^d}\right]\right)\left(\frac{\kappa}{\delta-\rho}\right)^{\frac{1}{\gamma}}
\end{eqnarray*}
and assumption $\left(U3\right)$ ensure that $\lim_{n\to \infty}u_{n,k}=0$, it then follows that $1-u_{n,k}>0$ for $n$ sufficiently large; therefore we can deduce that (see Assumption $\left(U4\right)i)$):
\begin{eqnarray*}
\mathbb{P}\left[\|s_n\|\geq w_n\quad \mbox{and}\quad v_n\left|\Psi_n\left(s_n\right)\right|\geq \delta\right] & \leq & \sum_{i=1}^n\mathbb{P}\left[\|X_k\|^{\beta}\geq w_n^{\beta}\left(1-u_{n,k}\right)^{\beta}\right]\\
&\leq & \sum_{i=1}^n\mathbb{E}\left(\|X_k\|^{\beta}\right) w_n^{-\beta}\left(1-u_{n,k}\right)^{-\beta}\\
&\leq & n \mathbb{E}\left(\|X_1\|^{\beta}\right) w_n^{-\beta}\max_{1\leq k\leq n}\left(1-u_{n,k}\right)^{-\beta}.
\end{eqnarray*}
Consequently, 
\begin{eqnarray*}
\lefteqn{\frac{\gamma_nv_n^2}{h_n^d}\log \mathbb{P}\left[\|s_n\|\geq w_n\quad \mbox{and}\quad v_n\left|\Psi_n\left(s_n\right)\right|\geq \delta\right]}\\
&&\leq  \frac{\gamma_nv_n^2}{h_n^d}\left[\log n+\log \mathbb{E}\left(\|X_1\|^{\beta}\right)-b\beta \frac{h_n^d}{\gamma_nv_n^2}-\beta\log \max_{1\leq k\leq n}\left(1-u_{n,k}\right)\right],
\end{eqnarray*}
and, thanks to assumptions $\left(U3\right)$, it follows that
\begin{eqnarray*}
\limsup_{n\to \infty}\frac{\gamma_nv_n^2}{h_n^d}\log \mathbb{P}\left[\|s_n\|\geq w_n\quad \mbox{and}\quad v_n\left|\Psi_n\left(s_n\right)\right|\geq \delta\right]\quad \leq \quad -b\beta,
\end{eqnarray*}
which concludes the proof of Lemma~\ref{lem7}.
\subsection{Proof of Proposition~\ref{pr:psiunif}}\label{proof:prop2}
Let us at first note that the lower bound
\begin{eqnarray}\label{eq:lowboun}
\liminf_{n\to \infty}\frac{\gamma_nv_n^2}{h_n^d}\log \mathbb{P}\left[\sup_{x\in U}v_n\left|\Psi_n\left(x\right)\right|\geq \delta\right]\geq -\tilde{g}_U\left(\delta\right)
\end{eqnarray}
follows from the application of Proposition~\ref{prop:LMDP} at a point $x_0\in \overline{U}$ such that $f\left(x_0\right)=\|f\|_{U,\infty}$.\\
In the case $U$ is bounded, Proposition~\ref{pr:psiunif} is thus a straightforward consequence of~(\ref{eq:lowboun}) and the first part of Lemma~\ref{lem5}. Let us now consider the case $U$ is unbounded.\\
Set $\delta>0$ and, for any $b>0$ set $w_n=\exp\left(b\frac{h_n^d}{\gamma_nv_n^2}\right)$. Since, by definition of $s_n$,
\begin{eqnarray*}
\lefteqn{\mathbb{P}\left[\sup_{x\in U}v_n\left|\Psi_n\left(x\right)\right|\geq \delta\right]}\\
&&\leq \mathbb{P}\left[\sup_{x\in U,\|x\|\leq w_n}v_n\left|\Psi_n\left(x\right)\right|\geq \delta\right]+\mathbb{P}\left[\|s_n\|\geq w_n \,\, \mbox{and}\,\, v_n\left|\Psi_n\left(x\right)\right|\geq \delta\right],
\end{eqnarray*} 
it follows from Lemmas \ref{lem5} and \ref{lem7} that 
\begin{eqnarray*}
\limsup_{n\to \infty}\frac{\gamma_nv_n^2}{h_n^d}\log \mathbb{P}\left[\sup_{x\in U}v_n\left|\Psi_n\left(x\right)\right|\geq \delta\right]\quad \leq \quad \max\left\{-b\beta;db-\tilde{g}_U\left(\delta\right)\right\}
\end{eqnarray*}
and consequently
\begin{eqnarray*}
\limsup_{n\to \infty}\frac{\gamma_nv_n^2}{h_n^d}\log \mathbb{P}\left[\sup_{x\in U}v_n\left|\Psi_n\left(x\right)\right|\geq \delta\right]\quad \leq \quad \inf_{b>0}\max\left\{-b\beta;db-\tilde{g}_U\left(\delta\right)\right\}.
\end{eqnarray*}
Since the infimum in the right-hand-side of the previous bound is achieved for $b=\tilde{g}_U\left(\delta\right)/\left(\beta+b\right)$ and equals $-\beta\tilde{g}_U/\left(\beta+d\right)$, we obtain the upper bound
\begin{eqnarray*}
\limsup_{n\to \infty}\frac{\gamma_nv_n^2}{h_n^d}\log \mathbb{P}\left[\sup_{x\in U}v_n\left|\Psi_n\left(x\right)\right|\geq \delta\right]\quad \leq \quad -\frac{\beta}{\beta+d}\tilde{g}_U\left(\delta\right)
\end{eqnarray*}
which concludes the proof of Proposition~\ref{pr:psiunif}.

\end{document}